\documentclass[12pt,letterpaper]{article}
\usepackage{amsfonts}
\usepackage{amssymb}
\usepackage{latexsym}  
\usepackage{pifont}  
\usepackage{marvosym}  

\input cyracc.def
\newfam\cyrfam

\font\twelvecyr=wncyr10 scaled 1200 
\def\cyr{\fam\cyrfam\twelvecyr\cyracc}

\renewcommand\AA{{\mathcal A}}
\newcommand\EE{{\mathcal E}}
\newcommand\FF{{\mathcal F}}
\newcommand\II{{\mathcal I}}
\newcommand\BB{{\mathcal B}}
\newcommand\MM{{\mathcal M}}
\newcommand\UU{{\mathcal U}}
\newcommand\HH{{\mathcal H}}
\newcommand\WW{{\mathcal W}}

\newcommand\RRR{{\mathbb R}}
\newcommand\NNN{{\mathbb N}}

\newcommand\QQQ{{\mathbb Q}}
\newcommand\TTT{{\mathbb T}}
\newcommand\CCC{{\mathbb C}}

\newcommand\KKKK{{\mathfrak K}}

\newcommand\cchi{{\raise 2 pt \hbox{$\chi$}}}

\newcommand\lparen{{\rm(}}
\newcommand\rparen{{\hskip 1pt \rm)}}

\newcommand\rest{\upharpoonright}     
\newcommand\res{\mathord {\upharpoonright}}  

\renewcommand\Re{\mathrm{Re}}   
\newcommand\diam{\mathrm{diam}}   
\newcommand\BIG{\mathrm{BIG}}   
\newcommand\supt{\mathrm{supt}}   
\newcommand\cl{\mathrm{cl}}   

\newcommand\No{\mathsf{n}}
\newcommand\Ea{\mathsf{e}}
\newcommand\So{\mathsf{s}}
\newcommand\We{\mathsf{w}}
\newcommand\Rt{\mathsf{r}}
\newcommand\Lft{\mathsf{l}}
\newcommand\Dgen{\mathsf{d}}  

\newcommand\Sh{\mbox{\cyr Sh}}  

\newcommand\sbl[1]{\langle#1\rangle}   

\newcommand\mosim{\mathord{\sim}}  

\newcommand\lea{\sqsubseteq}  
\newcommand\leac{\, \mbox{\rlap{$\sqsubseteq$}  %
   \raise 2.3pt \hbox{\scriptsize $\mathrm{c}$}}\; }

\newcommand\iv{^{-1}} 

\def\eop{{\Large \Coffeecup}}  

{\end{enumerate}}
\newenvironment{itemizz}{\begin{itemize}\setlength{\itemsep}{-1mm}} %
{\end{itemize}}                              

\newenvironment{itemizn}[1] 
{\begin{itemize} \setlength{\itemsep}{-1mm} %
\renewcommand\labelitemi{\ding{#1}}} %
{\end{itemize}}

\newtheorem{theorem}{Theorem}[section]
\newtheorem{definition}[theorem]{Definition}
\newtheorem{lemma}[theorem]{Lemma}
\newtheorem{corollary}[theorem]{Corollary}
\newtheorem{proposition}[theorem]{Proposition}
\newtheorem{example}[theorem]{Example}

\newenvironment{proof}{{\bf Proof.}}{\eop\medskip}
\newenvironment{proofof}[1]{\medskip \textbf{Proof of #1.}}{\eop\medskip}

\begin{document}

\title{
Complex Function Algebras and Removable Spaces}

\author{Joan E. Hart\footnote{University of Wisconsin, Oshkosh,
WI 54901, U.S.A.,
\ \ hartj@uwosh.edu}
\  and
Kenneth Kunen\footnote{University of Wisconsin,  Madison, WI  53706, U.S.A.,
\ \ kunen@math.wisc.edu}
\thanks{Both authors partially supported by NSF Grant DMS-0097881.}
}

\maketitle

\begin{abstract}
The compact Hausdorff space $X$ has the
\textit{Complex Stone-Weierstrass Property} (CSWP)
iff it satisfies the complex version of the Stone-Weier\-strass Theorem.
W. Rudin showed that all scattered spaces have the CSWP.
We describe some techniques for proving that certain non-scattered
spaces have the CSWP.  In particular, if $X$ is the product
of a compact ordered space and a compact scattered space,
then $X$ has the CSWP if and only if $X$ does not contain 
a copy of the Cantor set.
\end{abstract}

\section{Introduction} 
\label{sec-intro}
All topologies discussed in this paper are assumed to be Hausdorff.

\begin{definition}
\label{def-cswp}
If $X$ is compact, then
$C(X) = C(X,\CCC)$ is the algebra of continuous
complex-valued functions on $X$,
with the usual supremum norm.
$\AA \lea C(X)$ means that $\AA$ is a subalgebra of $C(X)$
which separates points and contains the constant functions.
$\AA \leac C(X)$ means that $\AA \lea C(X)$ and
$\AA$ is closed in $C(X)$.
$X$ has the \emph{Complex Stone-Weierstrass Property (CSWP)}
iff every $\AA \lea C(X)$ is dense in $C(X)$.
\end{definition}

Equivalently, $X$ has the CSWP iff every $\AA \leac C(X)$ equals $C(X)$.
Note that if we replaced ``$\CCC$'' by ``$\RRR$''
in Definition \ref{def-cswp}, the property would be true of all
compact $X$ by the Stone-Weierstrass Theorem.

By 1960, it was known that the CSWP is true
of all compact scattered spaces (Rudin \cite{RUD2}).
It was not known whether having the CSWP was equivalent to being scattered,
although two important examples were known of non-scattered
spaces which failed the CSWP, namely,
the Cantor set (Rudin \cite{RUD1}) and $\beta\NNN$
(Hoffman and  Singer \cite{HS}).
It was also well-known (and easy to see) that if $X$ fails the CSWP,
then so does every compact space containing $X$.
In particular, every compact space containing a Cantor subset
fails the CSWP.  It follows, as Rudin pointed out,
that for compact metric spaces, having the
CSWP is equivalent to being scattered.

These results are discussed in more detail in \cite{KU},
which shows that indeed there are non-scattered spaces with the CSWP.
Of course, none of these can contain a Cantor subset.

\begin{theorem}[\cite{KU}]
\label{thm-sep-lots}
If $X$ is a compact LOTS which does not contain a Cantor subset,
then $X$ has the CSWP.
\end{theorem}
As usual, a LOTS is a linearly ordered set with its order topology.
As a specific example,
the Aleksandrov double arrow space, which is not scattered, has the CSWP.
In this paper, we extend these
results to a much wider class of spaces.  Our results have,
as a special case:

\begin{theorem}
\label{thm-lots-scattered}
If $L$ is a compact LOTS which does not contain a Cantor subset,
and $S$ is a compact scattered space, then 
$L \times S$ has the CSWP.
\end{theorem}

Theorem \ref{thm-sep-lots} is the special case of Theorem
\ref{thm-lots-scattered} in which $S$ is a singleton.

We shall in fact proceed by a generalization of
the Cantor-Bendixson sequence.  The standard  Cantor-Bendixson sequence
is obtained
by removing isolated points:  Let $X'$ denote the set
of non-isolated points of $X$.  
Starting with a space $S$, we iterate this construction to 
obtain $S^{(0)} = S$, 
$S^{(\alpha + 1)} = (S^{(\alpha)})'$, and
$S^{(\gamma) } = \bigcap_{\alpha < \gamma}S^{(\alpha)}$
for limit ordinals $\gamma$. 
Then $S$ is \textit{scattered} iff some $S^{(\alpha)}$ is empty.

We shall define a class of 
pseudo-removable (PR) spaces.  There is a related class of PR-scattered
spaces where one obtains the empty set after repeatedly removing
open sets whose closure is PR (see Definitions \ref{def-PR} and
\ref{def-PR-scattered}).
We shall show (Theorem \ref{thm-PR-scattered-CSWP})
that every PR-scattered space has the CSWP.
The one-point space is PR, so that every scattered space is PR-scattered.
But also, we shall show (Lemma \ref{lemma-PR-scattered-lots})
that every compact LOTS which does not contain a Cantor
subset is PR, from which Theorem \ref{thm-lots-scattered}
will be obvious.

Section \ref{sec-basic} describes our basic techniques
and the outline of the main proofs, leaving some more
technical details to be verified in Sections \ref{sec-rem}
and \ref{sec-supt}.
There are still many open questions about the CSWP;
these are described in Section \ref{sec-conc}.

When we deal with scattered spaces in this paper,
it will be more useful to use the equivalent definition that every
non-empty set contains an isolated point.  Likewise, in the various
generalizations of scattered, we shall use a definition of form (1)
in the following simple observation:

\begin{proposition}
\label{prop-scattered}
Let $\KKKK$ be any class of compact spaces such that every
space homeomorphic to a
closed subspace of a member of $\KKKK$ is in $\KKKK$.  Then for
any compact $X$, the following are equivalent:
\begin{itemizz}
\item[1.] For all non-empty closed $F \subseteq X$,
there is a non-empty $U \subseteq F$
which is (relatively) open in $F$ such that $\overline U \in \KKKK$.
\item[2.] If one defines $Z^*$ to be
$Z \setminus \bigcup\{U \subseteq Z: U \mbox{\rm\  is open in $Z$ and }
\overline U \in \KKKK\}$, and then defines
$X^{[0]} = X$, 
$X^{[\alpha + 1]} = (X^{[\alpha]})^*$, and
$X^{[\gamma] } = \bigcap_{\alpha < \gamma}X^{[\alpha]}$
for limit ordinals $\gamma$,
then some $X^{[\alpha]}$ is empty.
\end{itemizz} 
\end{proposition}
\begin{proof}
For $(2) \to (1)$, consider the largest $\alpha$ such that 
$F \subseteq X^{[\alpha] }$.
\end{proof}

\begin{definition}
\label{def-scattered}
If $\KKKK$ is as in Proposition \ref{prop-scattered} and $X$ is compact,
then $X$ is \emph{scattered for} $\KKKK$ iff $X$ satisfies
$(1)$ \lparen or $(2)$\rparen.
\end{definition}

So, scattered spaces in the usual sense are scattered for
the class of 1-point spaces.

\section{Basic Techniques}
\label{sec-basic}
These techniques involve \textit{idempotents}, \textit{measures},
and \textit{removable spaces}.

As usual, $f \in C(X)$ is called an \textit{idempotent} iff
$f^2 = f$; equivalently, $f = \cchi_H$, where $H$ is a clopen
subset of $X$ (i.e., $H$ is both closed and open).  Thus,
$X$ is connected iff the only idempotents in $C(X)$ are the
two \textit{trivial} ones (the constant $0$ and the constant $1$ functions).

\begin{definition}
The compact $X$ has the \emph{NTIP} iff every $\AA \leac C(X)$
contains a non-trivial idempotent.
\end{definition}

A connected space cannot have the NTIP, but if $X$ is not connected,
the CSWP implies the NTIP.  Rudin \cite{RUD1} proved that the Cantor
set fails the NTIP (and hence fails the CSWP).
By \cite{KU}, Lemma 3.5 (or see Proposition \ref{prop-de-branges}):

\begin{lemma}
\label{lemma-ntip-cswp}
If $X$ is compact and every perfect subset of $X$ has the NTIP,
then $X$ has the CSWP.
\end{lemma}

This simplifies the task of proving that a space has the CSWP,
since one need only produce a non-trivial idempotent,
although one must then deal with arbitrary perfect subsets of $X$.
One can produce idempotents in some $\AA \leac C(X)$ by
using Runge's Theorem on polynomial approximations (see \cite{RUD3}\S13).
We quote here only the special case we need:

\begin{theorem}[Runge]
\label{thm-runge}
If $K_1 , \ldots , K_n$ are disjoint compact convex subsets of $\CCC$,
$\varepsilon > 0$, and
$w_1, \ldots, w_n \in \CCC$,
then is a complex polynomial $P(z)$ such that 
$|P(z) - w_\ell| < \varepsilon$
for each $\ell = 1, \ldots, n$ and each $z \in K_\ell$.
\end{theorem}

Then, as in \cite{RUD1} or \cite{HS}, composing functions with
polynomials yields:

\begin{lemma}
\label{lemma-disc}
If $X$ is compact, $\AA \leac C(X)$, and $\Re(h(X))$ is not
connected for some $h \in \AA$, then $\AA$ contains a non-trivial
idempotent.
\end{lemma}

As usual, $\Re: \CCC \to \RRR$ denotes projection onto the real axis.
Lemmas \ref{lemma-ntip-cswp} and \ref{lemma-disc} were used in \cite{KU}
to prove Theorem \ref{thm-sep-lots} in the case that $X$ is separable.

Note that any $X$ satisfying the
hypotheses to Lemma \ref{lemma-ntip-cswp} is totally disconnected.
However, using \textit{measures}, we can extend our results
to apply to many connected spaces.
A simple use of measures is contained in \cite{KU};
Lemma \ref{lemma-supt-cswp} below is Lemma 5.2 of \cite{KU}.

\begin{definition}
If $\mu$ is a regular complex Borel measure on the compact space $X$,
then $|\mu|$ denotes its total variation, and 
$\supt(\mu) = \supt(|\mu|)$ denotes its (closed) support;
that is,
$\supt(\mu) = X \setminus \bigcup \{U\subseteq X :\; $U$ \mbox{ is open }
\ \ \& \ \ |\mu|(U) = 0\}$.
\end{definition}

\begin{lemma}
\label{lemma-supt-cswp}
Assume that $X$ is compact and that $\supt(\mu)$ has the CSWP
for all regular Borel measure $\mu$.
Then $X$ has the CSWP.
\end{lemma}

In fact, one can derive Lemmas
\ref{lemma-supt-cswp} and \ref{lemma-ntip-cswp} together
by the method of  Branges \cite{BRA};
see Proposition \ref{prop-de-branges}.
In \cite{KU}, Lemma \ref{lemma-supt-cswp}
was used to derive Theorem \ref{thm-sep-lots} 
from the separable case of it, by applying:

\begin{lemma}
\label{lemma-lots-supt}
Assume that $X$ is a compact LOTS and that $\mu$ is a regular Borel measure.
Then $\supt(\mu)$ is separable.
\end{lemma}

In Section \ref{sec-supt},
we shall continue to use the notion of ``support'' to 
reduce the problem of the CSWP for a ``big'' space $X$ to
that of a ``small'' subspace.
Usually, the ``small'' subspace will be totally disconnected,
and arguments involving idempotents will apply to it,
whereas $X$ itself may be connected.
For example, it is easy to construct
a connected compact LOTS which
does not contain a Cantor subset (see Example \ref{ex-PR-lots}),
but every separable subspace of such a LOTS must
be totally disconnected.

Theorem \ref{thm-lots-scattered} will be proved using our 
notion of \textit{removable spaces}.
The definition of ``removable'' is given in terms of the \v Silov boundary:

\begin{definition}
\label{def-silov}
Assume that $\AA\leac C(X)$.  Let $H$ be a closed subset of $X$.  Then:
\begin{itemizn}{"2B}
\item
$\|f\|_H = \sup\{|f(x)| : x \in H\}$.
\item
$H$ is a \emph{boundary} for $\AA$ iff
$\|f\|_H = \|f\|$ for all $f\in \AA$.
\item
$\Sh(\AA)$ denotes the  \v Silov \emph{boundary};
this is the smallest closed set which is a boundary for $\AA$.
\end{itemizn}
\end{definition}
The existence of a smallest boundary, $\Sh(\AA)$,
is due to \v Silov; see \cite{AW, GA, HOR, RUD4}.
Note that $\Sh(\AA)$ is always non-empty, and cannot be finite
unless $X$ is finite.

\begin{definition}
Given $\AA \lea C(X)$ and a closed subset $H \subseteq X$,
let $\AA\res H = \{f \res H : f \in \AA\} \lea C(H)$.
$\overline{\AA\res H} \leac C(H) $ denotes the closure in the uniform topology.
\end{definition}

Observe that $\AA \leac C(X)$ does not in general imply
that $\AA\res H  \leac C(H)$.  For example, if 
$X \subset \CCC$ is the closed unit disc and $\AA$ is the usual disc algebra
(see \cite{GA, RUD4}),
then $\AA \res [0,1]$ is dense in and not closed in $C([0,1])$.
The following easy lemma is proved in  \cite{KU}:

\begin{lemma}
\label{lemma-restrict-shilov}
Suppose that $\AA \leac C(X)$ and $\Sh(\AA) \subseteq H$, where $H$ is closed.
Then $\AA\res H  \leac C(H)$.  Also, 
$\AA\res H  = C(H)$ iff
$\AA = C(X)$ \lparen in which case $\Sh(\AA) = X$\rparen.
If all idempotents of $\AA$ are trivial, then the same is true
of $\AA\res H$.
\end{lemma}

\begin{definition}
\label{def-remove}
A compact space $H$ is \emph{removable} iff for all $X,U,\AA$, if:
\begin{itemizn}{"2B}
\item $X$ is compact,
\item $U \subsetneqq X$ and $U$ is open,
\item $\overline U$ is homeomorphic to a subspace of $H$, and
\item $\AA \leac C(X)$ and all idempotents of $\AA$ are trivial,
\end{itemizn}
then $\Sh(\AA) \subseteq X\backslash U$.
\end{definition}

We remark that we disallow $U = X$ only because if $X$ is a singleton,
$U = X$ and $\AA = C(X)$, then we would have a contradiction.
If $U = X$ and $X$ is not a singleton, then the conclusion,
$\Sh(\AA) = \emptyset$, is still contradictory, but so
are the hypotheses, since $X$ would have the  NTIP by
Lemma \ref{lemma-remove-basic} below.
As stated, the hypotheses are non-vacuous, in the sense that
if $H$ is any compact space, then there are always 
$X,U,\AA$ satisfying the hypotheses of Definition \ref{def-remove}
with $U = \overline U$ homeomorphic to $H$; see
Proposition \ref{prop-compact-copy}.  

\begin{lemma}
\label{lemma-remove-basic}
Let $H$ be removable.
\begin{itemizz}
\item[1.] If $K \subseteq H$ is closed, then $K$ is removable.
\item[2.] If $|H| > 1$, then $H$ has the NTIP.
\item[3.] $H$ has the CSWP.
\item[4.] $H$ is totally disconnected.
\end{itemizz}
\end{lemma}
\begin{proof}
(1) is immediate from the definition.  For (2), suppose that we had
$\AA \leac C(H)$ with all idempotents trivial.  Let $H = U \cup V$,
where $U,V$ are proper open subsets.  Applying (1), 
$\overline U$ and $\overline V$ are removable, so 
$\Sh(\AA) \subseteq X\backslash U$ and
$\Sh(\AA) \subseteq X\backslash V$, so $\Sh(\AA) = \emptyset$,
a contradiction.

Now (3) is immediate from (2), (1), and Lemma \ref{lemma-ntip-cswp}.
(4) is also immediate from (2) and (1), since no
perfect subset of $H$ is connected.
\end{proof}

Of course, for this notion to be useful, we need to show that 
there are some removable spaces.  We begin with:

\begin{lemma}
\label{remov-pt}
The one-point space is removable.
\end{lemma}
\begin{proof}
Let $X,U,\AA$ be as in Definition \ref{def-remove}, with
$\overline U$ a singleton, $\{p\}$, so that $p$ is isolated in $X$.
We need to show that $X \backslash \{p\}$ is a boundary.
If not, then there is an $h\in\AA$ with $\|h\|_{ X \backslash \{p\} } = 1$
but $|h(p)| > 1$.  We may assume that
$h(p)\in\RRR$; but then 
$\Re(h(X))$ is not connected, a contradiction by Lemma \ref{lemma-disc}.
\end{proof}

We can produce more removable spaces via a
generalized Cantor-Bendixson analysis:

\begin{definition}
The compact $H$ is \emph{R-scattered} iff $H$ is scattered for
the class of removable spaces
\lparen see Definition \ref{def-scattered}\rparen.
\end{definition}

Note that every closed subspace of an R-scattered space is R-scattered.
By Lemma \ref{remov-pt}, every scattered space is R-scattered,
and hence removable by:

\begin{lemma}
\label{lemma-R-scat-remov}
$H$ is R-scattered iff $H$ is removable.
\end{lemma}
\begin{proof}
For the non-trivial direction,
let $X,U,\AA$ be as in Definition \ref{def-remove}.
Then $\overline U$ is R-scattered.  Let $K = \Sh(\AA)$.  We need
to show that $U \cap K =  \emptyset$, so assume that 
$U \cap K \ne  \emptyset$.
Note that $X$ and  $K$ must be infinite.

By Lemma \ref{lemma-restrict-shilov}, $\AA\res K \leac C(K)$
and $\AA\res K$ has no non-trivial idempotents.
Also note that $\Sh(\AA\res K) = K$.

Let $W = U \cap K$.  Then $W$ is relatively open in $K$, and $W$ is
non-empty.  $\overline W \subseteq \overline U$ and
$\overline U$ is R-scattered, so choose $V \subseteq \overline W$
with $V$ non-empty, $V$ open in $\overline W$, and $\overline V$ removable.
We may assume also that $V \subseteq  W$ (otherwise,
replace $V$ by $V \cap W$), so that $V$ is open in $K$.
Now $V \ne K$ (since otherwise $K$ would have the NTIP by 
Lemma  \ref{lemma-remove-basic}),
so $\Sh(\AA \res K) \subseteq K \backslash V$, a contradiction.
\end{proof}

It follows that every scattered space is removable, and thus
has the CSWP, which was already known from Rudin \cite{RUD2};
in fact, all we have done is to redo the argument in \cite{RUD2},
using somewhat more complicated terminology.
For these methods to produce anything
new, we need to produce some non-scattered removable spaces,
which we do in Section \ref{sec-rem}.  In particular, we shall
prove there:

\begin{lemma}
\label{lemma-R-scattered-lots}
Every compact separable LOTS which does not contain a Cantor
subset is removable.
\end{lemma}

\begin{proofof}{Theorem \ref{thm-lots-scattered}}
Fix $X = L \times S$ as in Theorem \ref{thm-lots-scattered}.

Applying Lemma \ref{lemma-supt-cswp}, 
it is sufficient to fix a regular Borel measure $\mu$ on $X$ and
prove that $Y := \supt(\mu)$ has the CSWP.
For this, it is sufficient to show that $Y$ is R-scattered
(and hence removable by Lemma \ref{lemma-R-scat-remov}).

Let $\pi_L: L \times S \to L$ be projection, and let
$\mu\, \pi_L\iv$ be the induced measure on $L$.
Let $L' = \supt(\mu\, \pi_L\iv)$.  Then $L'$ is separable
(by Lemma \ref{lemma-lots-supt}), and 
$Y \subseteq L' \times S$.

To prove that $Y$ is R-scattered, fix a non-empty closed $F \subseteq Y$.
Since $S$ is scattered, $F$ contains a non-empty clopen subset of
the form $L'' \times \{p\}$, where $p$ is isolated in $\pi_S(F)$
and $L'' \subseteq L'$.  This  $L''$ is removable by
Lemma \ref{lemma-R-scattered-lots} (a subspace of a separable LOTS
is separable; see \cite{LB}).
\end{proofof}

We can generalize Theorem \ref{thm-lots-scattered} to
Theorem \ref{thm-PR-scattered-CSWP} below, using the following:

\begin{definition}
\label{def-PR}
The compact $H$ is \emph{pseudo-removable (PR)} iff
$\supt(\mu)$ is removable whenever $\mu$ is a regular
Borel measure on $H$.
\end{definition}

\begin{definition}
\label{def-PR-scattered}
The compact $X$ is \emph{PR-scattered} iff $H$ is scattered for
the class of PR spaces
\lparen see Definition \ref{def-scattered}\rparen.
\end{definition}

For example, applying Lemmas  \ref{lemma-lots-supt}
and \ref{lemma-R-scattered-lots}, we get:

\begin{lemma}
\label{lemma-PR-scattered-lots}
Every compact LOTS which does not contain a Cantor subset is PR.
\end{lemma}

We shall prove:

\begin{theorem}
\label{thm-PR-scattered-CSWP}
Every compact PR-scattered space has the CSWP.
\end{theorem}

Given this, the proof of Theorem \ref{thm-lots-scattered} 
is trivial, since it is immediate from 
Lemma \ref{lemma-PR-scattered-lots} that
$L \times S$ is PR-scattered.
Now, our proof above of Theorem \ref{thm-lots-scattered}
actually showed that $L \times S$ is PR,
but there are PR-scattered spaces which are not PR
(see Example \ref{ex-PR-lots}).  For these spaces,
the proof we gave of Theorem \ref{thm-lots-scattered} will not
work, but we need to use Theorem \ref{thm-PR-scattered-CSWP}
instead.

Of course, we still need to prove Lemma \ref{lemma-R-scattered-lots}
(in Section \ref{sec-rem}) and 
Theorem \ref{thm-PR-scattered-CSWP} (in Section \ref{sec-supt}).
Theorem \ref{thm-PR-scattered-CSWP} makes use of the following easy remark:

\begin{lemma}
\label{lemma-reduce-PR-scattered}
There is no PR-scattered compact space $F$
with more than one point,
with a regular Borel measure $\nu$ and an
$\AA \leac C(F)$ such that:
\begin{itemizz}
\item[1.] $F = \supt(\nu)$.
\item[2.] All idempotents of $\AA$ are trivial.
\item[3.] $F = \Sh(\AA)$.
\end{itemizz}
\end{lemma}
\begin{proof}
Let $U \subseteq F$ be non-empty and
open, such that $\overline U$ is PR.
$\overline U$ is actually removable, since it
is the support of the measure $\nu$ restricted to $\overline U$. 
Then $\overline U$ cannot equal $F$ (otherwise, $F$ would
have the NTIP by Lemma \ref{lemma-remove-basic}).
It follows that $\Sh(\AA) \subseteq F \backslash U$,
contradicting (3).
\end{proof}

The proof of Theorem \ref{thm-PR-scattered-CSWP}
takes a PR-scattered space $X$ and an $\AA \leac C(X)$
such that $\AA \ne C(X)$,
and produces an $F \subseteq X$ and a $\nu$ for which (1)(2)(3) hold for
$\overline{\AA \res F} \leac C(F)$,
which is contradictory. Now the methods discussed above can
easily get $F,\nu,\overline{\AA \res F}$ to satisfy (2)(3),
and it is also easy to get (1)(2) (see Proposition \ref{prop-de-branges}),
but obtaining (1)(2)(3) seems to require a different idea;
we shall obtain $\nu$ as a measure which represents an
element of the maximal ideal space of $\AA$.
The details are described in Section \ref{sec-supt}.

Sections \ref{sec-rem} and \ref{sec-supt} can be read 
independently of each other.

\section{Some Removable Spaces}
\label{sec-rem}
We describe a general technique for proving that
certain spaces are removable.  We begin with the observation,
following Tychonov, that any subset $\EE \subseteq C(X)$ 
defines a map from $X$ into $\CCC^\EE$:

\begin{definition}
\label{def-quot}
Let $X$ be compact and $\EE \subseteq C(X)$.
Define $x \sim_\EE y$ (or, $x \sim y$) iff $f(x) = f(y)$ for all $f\in\EE$.
Let $[x] = [x]_\EE = \{y : x \sim y\}$, let
$X/\mosim = \{[x] : x \in X\}$, and let $\pi = \pi_\EE$ be the natural
map from $X$ onto $X/\mosim$.  A subset $U \subseteq X$ is $\EE$--\emph{open}
iff $U$ is open and $U = \pi\iv \pi(U)$,
and $H \subseteq X$ is $\EE$--\emph{closed} iff $X \backslash H$ is
$\EE$--open.
For $U \subseteq X$, let $\cl_\EE(U) = \pi\iv( \cl(\pi(U)))$.
\end{definition}

With the usual quotient topology, $X/\mosim$ is a compact
space and can be identified with the image of $X$ under the
evaluation map from $X$ into $\CCC^\EE$.
Each $[x]$ is $\EE$--closed.
If $\EE$ is countable, then $X/\mosim$ will be second countable
(equivalently, metrizable).
$\cl_\EE(U)$ is the smallest $\EE$--closed set containing $U$.
Note that even for
$\EE$--open $U$, $\cl_\EE(U)$ might properly contain $\overline U$.

\begin{definition}
Let $X$ be compact, fix $\EE \subseteq C(X)$ and
$f \in C(X)$, and fix a real $c > 0$.
Then $\BIG_\EE(X,f , c) = \bigcup\{[x] \in X/\mosim: \diam(f([x])) \ge c\}$.
\end{definition}

Here, ``$\diam$'' refers to the
usual notion of the diameter of a subset of $\CCC$.
Note that $\BIG_\EE(X,f , c) = \emptyset$ whenever $f\in\EE$.

\begin{lemma}
\label{lemma-big-closed}
Each $\BIG_\EE(X,f , c)$ is closed in $X$.
\end{lemma}
\begin{proof}
Let $B = \{(x,y,z) \in X^3: x \sim y \sim z \ \&\ |f(y) - f(z)| \ge c\}
\subseteq X^3$.
Then $B$ is closed and $\BIG_\EE(X,f , c)$ is the projection
of $B$ onto the first coordinate.
\end{proof}

Using these notions, we can define a class of spaces which
are ``close'' to being metrizable:

\begin{definition}
\label{def-nice}
Let $X$ be compact and $\EE \subseteq C(X)$.
$\EE$ is \emph{nice} iff
\begin{itemizz}
\item[1.]
$[x]$ is scattered for all $[x] \in X$.
\item [2.]
For all $f\in C(X)$,
$f([x])$ is a singleton for all but at most countably
many equivalence classes $[x] \in X/\mosim$.
\end{itemizz}
$X$ is \emph{nice} iff there is \emph{countable} nice $\EE \subseteq C(X)$.
\end{definition}

\begin{lemma}
\label{lemma-big-scattered}
If $\EE \subseteq C(X)$ is nice, $c > 0$, and $f \in C(X)$, 
then $\BIG_\EE(X,f , c)$ is scattered.
\end{lemma}
\begin{proof}
It is compact by Lemma \ref{lemma-big-closed}, and
a countable union of scattered subspaces.
\end{proof}

We shall show (Lemma \ref{lemma-nice-sh})
that if $X$ is nice and does not contain a Cantor subset, then $X$ is removable.
Some examples of nice $X$:  If $X$ is scattered then $X$ is nice,
taking $\EE$ to contain only constant functions
(so there is only one equivalence class).
If $X$ is second countable, then $X$ is nice, since there is a
countable $\EE\subset C(X)$ which separates points, so that
each $[x]$ is a singleton.  Of course, any second countable $X$ is also
scattered if it does not contain a Cantor subset.  A more useful
example is:

\begin{lemma}
\label{lemma-sep-LOTS-nice}
Every separable compact LOTS is nice.
\end{lemma}
\begin{proof}
Let $D \subseteq X$ be dense in $X$ and countable.
Assume that $X$ and $D$ are infinite.
List $D$ as $\{d_n : n\in\NNN\}$.
We may assume that $d_0$ is the first element of $X$ and $d_1$ 
is the last element of $X$.
Define $h: D \to [0,1]$ so that $h(d_0) = 0$, $h(d_1) = 1$, and
$h(d_n) = (h(d_i) + h(d_j))/2$ when $n \ge 2$, where
$d_i$ is the largest element in $\{d_\ell : \ell < n \ \&\ d_\ell < d_n\}$, and
$d_j$ is the smallest element in $\{d_\ell : \ell < n \ \&\ d_\ell > d_n\}$.

This $h$ extends to a continuous map from $X$ into $[0,1]$,
defined so that for $x \in X \backslash D$, 
$\; h(x) = \sup\{h(d_\ell) : d_\ell < x\} = \inf\{h(d_\ell) : d_\ell > x\}$.
Let $\EE = \{h\}$. 
(1) of Definition \ref{def-nice} holds because
each $[x]$ has cardinality one or two.  To prove (2), note that
for any $f\in C(X)$ and
each $c >0$,
$\BIG_\EE(X,f , c)$ must be a \textit{finite} union of two-element classes,
since if it were infinite, it would have a limit point,
which would contradict continuity of $f$.
\end{proof}

A class of nice spaces which is not related to ordered spaces
or to scattered spaces is described in Example \ref{ex-nice-not-ordered}.

The definitions of ``nice'' and ``$\BIG$'' refer only to $\sim_\EE$,
not $\EE$, so it is reasonable to introduce:

\begin{definition}
If $\EE,\FF \subseteq C(X)$, say that
$\EE \preccurlyeq \FF$ iff $\sim_\FF$ is finer than $\sim_\EE$
(that is, every $\sim_\EE$ class is a union of $\sim_\FF$ classes).
Say $\EE \approx \FF$ iff $\sim_\FF$ and $\sim_\EE$ are the same.
\end{definition}

Note that $\EE \subseteq \FF \to \EE \preccurlyeq \FF$, and
$\EE \approx \FF \leftrightarrow \EE \preccurlyeq \FF \preccurlyeq \EE$.

\begin{lemma}
\label{lemma-refine-nice}
If $\EE \preccurlyeq \FF$,  $\EE$ is nice, and $\FF$ is countable,
then $\FF$ is nice.
\end{lemma}
\begin{proof}
In Definition \ref{def-nice}, clause (1) for $\FF$ is obvious.
To verify clause (2), note that every $\EE$ class, $[x]_\EE$
is a countable union of $\FF$ classes because
$([x]_\EE)/\mosim_\FF$ is scattered and second countable,
and hence countable.
\end{proof}

Next, we mention a few closure properties of the class of nice space:

\begin{lemma}
If $X$ is nice, then every closed subspace of $X$ is nice.
\end{lemma}

\begin{lemma}
\label{lemma-union-nice}
If $X$ is compact and $X = \bigcup_{n\in\omega} H_n$, where each $H_n$
is nice and is a closed $G_\delta$ in $X$, then $X$ is nice.
\end{lemma}
\begin{proof}
For each $n$, choose a countable nice $\EE_n \subseteq C(H_n)$.
By the Tietze Extension Theorem, we may assume that
$\EE_n = \FF_n \res H_n$, where $\FF_n \subseteq C(X)$.
Since $H_n$ is a $G_\delta$,
we may also assume that each $H_n$ is $\FF_n$--closed
(adding another function to $\FF_n$ if necessary).
Then $\bigcup_n  \FF_n \subseteq C(X)$ is countable and nice.
\end{proof}

\begin{corollary}
\label{cor-sum-nice}
A finite disjoint sum of nice spaces is nice.
\end{corollary}

\begin{corollary}
If $X$ is nice and
$Y$ is countable and compact, then $X \times Y$ is nice.
\end{corollary}

The product of two nice spaces is not in general nice;
for example, the square of the double arrow space
is not nice by Proposition \ref{prop-bad-product}.

Next, we relate ``nice'' to function algebras.

\begin{definition}
If $X$ is compact and $\EE \subseteq C(X)$, then
$\sbl{\EE}$ is the intersection of all closed subalgebras
of $C(X)$ which
contain $\EE$ and the constant functions.
\end{definition}

So, $\sbl{\EE}$ is, by definition, a closed subalgebra which contains
the constant functions, but it will not separate points unless $\EE$ does:

\begin{lemma}
\label{lemma-expandE}
Let $X$ be compact,  $\EE \subseteq C(X)$, and 
$\EE^*$ the set of all complex conjugates of functions in $\EE$.
Let $\AA = \sbl{\EE}$ and $\AA' = \sbl{\EE \cup \EE^*}$.
Then $\EE \approx \AA \approx \AA'$.
\end{lemma}

\begin{lemma}
\label{lemma-move-E}
Say $\EE \subseteq C(X)$ is countable and $\AA \lea C(X)$.
Then there is a countable $\FF \subseteq \AA$ such $\EE \preccurlyeq \FF$.
\end{lemma}
\begin{proof}
By the Stone-Weierstrass Theorem, we may choose a
countable $\FF \subseteq \AA$ such that 
$\EE \subseteq \sbl{\FF \cup \FF^*}$.
Now, apply Lemma \ref{lemma-expandE}.
\end{proof}

Applying Lemma \ref{lemma-move-E} and \ref{lemma-refine-nice},
we get

\begin{lemma}
\label{lemma-inside-A}
If $X$ is nice and $\AA \lea C(X)$, then there is a
countable nice $\FF \subseteq \AA$.
\end{lemma}

In order to relate ``nice'' to ``removable'', we need
to know that the nice $\EE$ may be taken to have the
following additional property:

\begin{definition}
\label{def-adequate}
Suppose that $\EE \subseteq \AA \lea C(X)$.  Then
$\EE$ is \emph{adequate} for $\AA$ iff
whenever we are  given a finite $m$, an $f \in \AA$,
and sets $H_i,W_i$ for $i < m$
satisfying:
\begin{itemizz}
\item[1.] The $H_i$ are $\EE$--closed subsets of $X$,
\item[2.] The $W_i$ are open subsets of $\CCC$, and
\item[3.] $f(H_i) \subseteq W_i$ for each $i< m$,
\end{itemizz}
then there is a $g\in\EE$ such that $g(H_i) \subseteq W_i$ for each $i< m$.
\end{definition}

\begin{lemma}
\label{lemma-get-adequate}
If $\EE \subseteq \AA \lea C(X)$ and $\EE$ is countable,
then there is a countable $\FF$ such that
$\EE \subseteq \FF \subseteq \AA$ and 
$\FF$ is adequate for $\AA$.
\end{lemma}
\begin{proof}
We first reduce ``adequate'' to a countable number 
of instances.  Let $\WW$ be a countable base for the topology
of $\CCC$.  Assume that $\WW$ is closed under finite unions.
Let $\UU_\EE$ be a countable base for the space $X/\mosim$, and also assume
that $\UU_\EE$ is closed under finite unions.
Let $\pi = \pi_\EE$ (see Definition \ref{def-quot}), and let
$\HH_\EE =\{X  \setminus  \pi\iv(U) : U \in \UU_\EE\}$.
Then all sets in $\HH_\EE$ are $\EE$--closed.

Note that in order for $\EE$ to be adequate for $\AA$,
it is sufficient to verify Definition \ref{def-adequate}
with the $H_i \in \HH_\EE$ and the
$W_i \in \WW$:  To see this, suppose we started with 
arbitrary $H_i$ and $W_i$.
Each $H_i$ is the intersection of the sets in $\HH_\EE$ which contain it,
so by compactness, we can find $H'_i \supseteq H_i$ in $\HH_\EE$
with $f(H'_i) \subseteq W_i$.  Again by compactness,
we can find $W'_i \subseteq W_i$ in $\WW$
with $f(H'_i) \subseteq W'_i$.
Then if we get $g \in \EE$ 
with each $g(H'_i) \subseteq W'_i$, we will also have
$g(H_i) \subseteq W_i$.

Now, starting from $\EE$, we get
$\EE = \EE_0 \subseteq \EE_1 \subseteq \EE_2 \subseteq \cdots$,
where each $\EE_n$ is countable.  Given $\EE_n$, we obtain
$\EE_{n+1}$ so that whenever  $m,f $ and the
$H_i,W_i$ satisfy (1)(2)(3) of Definition \ref{def-adequate}
with the $H_i \in \HH_{\EE_n}$ and the $W_i \in \WW$,
there is a $g\in\EE_{n+1}$ such that
$g(H_i) \subseteq W_i$ for each $i< m$.
Let $\FF = \bigcup_n \EE_n$.
\end{proof}

Finally, we need an easy consequence of Runge's Theorem
\ref{thm-runge}:

\begin{lemma}
\label{lemma-runge}
Suppose that $\AA \lea C(X)$, $V \subseteq X$, $p,q \in X$,
$p \ne q$,
and $\AA$ contains a function $k$ with
$|k(p)|, |k(q)| > \sup\{|k(x)| : x \in V\}$.
Fix $a,b,c \in \CCC$ and $\varepsilon > 0$.
Then $\AA$ contains a function $f$ with
$f(p) = a$, $f(q) = b$, and $f(V) \subseteq B(c;\varepsilon)$.
\end{lemma}
\begin{proof}
Since $\AA$ separates points, we may assume, by adding a small function
to $k$, that $k(p) \ne k(q)$.  Then, whenever $\delta > 0$, we may
apply Theorem \ref{thm-runge} to find a polynomial $P$
such that, setting $h = P\circ k$, we have
$h(p) \in B(a; \delta)$,
$h(q) \in B(b; \delta)$, and $h(V) \subseteq B(c;\delta)$.
The result now follows by choosing $\delta$ small enough
and composing $h$ with a linear polynomial.
\end{proof}

\begin{lemma}
\label{lemma-nice-sh}
Assume that $H$ is nice and does not contain a Cantor subset.
Then $H$ is removable.
\end{lemma}
\begin{proof}
Fix  $X,U,\AA$ satisfying the conditions in Definition \ref{def-remove}.
Then  $\overline U$ is nice and $X\backslash U \ne \emptyset$, and
we need to prove that $X\backslash U$ is a boundary.  Assume
that it is not a boundary, and we shall derive a contradiction
by producing a non-trivial idempotent in $\AA$.
Fix a $k \in \AA$
such that $\|k\| > 1$ but $\|k\|_{X\backslash U} \le 1$.
By replacing $k$ with a power of $k$, we assume also that
$\|k\| > 3$.  Multiplying $k$ by some $e^{i\theta}$, we can
assume that some $k(x)$ is a real number in $(3, \infty)$.

Let $T = \{x\in X : |k(x)| \ge 2\}$ and
$V = X \backslash T= \{x\in X : |k(x)| < 2 \}$.
Then $T \subseteq U$, so $T$ is nice.
Applying Lemma \ref{lemma-inside-A} to $\AA \res T$,
let $\EE_0 \subset \AA$ be countable,
with $\EE_0 \res T \subset C(T)$ nice.
Then by Lemma \ref{lemma-get-adequate},
let $\EE$ be countable with $\EE_0 \cup \{k\} \subseteq \EE \subset \AA$ and
$\EE$ adequate for $\AA$.  Note that $\EE \res T \subset C(T)$
is nice by Lemma \ref{lemma-refine-nice}.
In the following, $[x]$ and $\sim$ always means $[x]_\EE$ and $\sim_\EE$.
Since $k\in\EE$, each $[x]$ is a subset of either $V$ or $T$.

$\{x \in X : |k(x)| \ge 3\}$ is not scattered
(since  $\Re(k(X))$ is connected by Lemma \ref{lemma-disc}), so
it cannot be second countable
(since it does not contain a Cantor subset),
so fix distinct $p,q\in \{x \in X : |k(x)| \ge 3\}$ with $p \sim q$.  Then, 
by Lemma \ref{lemma-runge},  fix
$f \in \AA$ with $f(V) \subseteq B(0; 1/2)$ and $f(p) = 2$
and $f(q) = -2$.

Let $H = \BIG_\sim(T,f\res T, 1/2) \subseteq T$.
$\; H$ is scattered by Lemma \ref{lemma-big-scattered}.
It follows that $f(H) - f(H) = \{f(x) - f(y) : x,y \in H\}$
is a continuous image of the scattered space $H \times H$,
and is thus scattered in $\CCC$, so fix
a $b \in (1, 2)$ with $ b \ne \Re(f(H) - f(H))$.

Then, $\pm b \ne \Re(f([x]) - f([x])$ for each $x \in X$:
If $x \in H$, this follows by the choice of $b$.
If $x\in T \backslash H$, then $\diam(f([x])) < 1/2$,
so that  $f([x]) - f([x]) \subseteq B(0; 1)$.
If $x \in V$, then
$f([x]) \subseteq f(V) \subseteq B(0; 1/2)$.

Next, for each $x\in X$, there are $\EE$--open $U = U_x$ containing
$[x]$ and open $W = W_x \subseteq \CCC$ such that
$f(\cl_\EE(U)) \subseteq W$
and $\pm b \notin \Re(W - W)$.
To do this, choose $\EE$-open $U_n$  so that
$U_0 \supseteq \cl_\EE( U_1)  \supseteq U_1 \supseteq \cl_\EE( U_2) \cdots $
and $\bigcap_n U_n = [x]$.   Since $\pm b \notin \Re(f([x]) - f([x]))$,
we can apply compactness to choose $n$
so that $\pm b \notin \Re(f(\cl_\EE (U_n)) - f(\cl_\EE (U_n))$.
Then let $U_x$ be that $U_n$, and let $W$ be any open superset of 
$f(\cl_\EE (U_n))$ such that $\pm b \notin \Re(W - W)$.

Then, by compactness, we have a finite $m$ and $U^i = U_{x_i}$ for $i < m$
such that the $U^i$ cover $X$.
We may assume that $x_0\in V$ and 
$U^0 = V$, and $W_{x_0} = B(0; 1/2)$.
Let $W^i = W_{x_i}$.
$f(\cl_\EE(U^i)) \subseteq W^i$ and $\pm b \notin \Re(W^i - W^i)$.
Since $\EE$ is adequate,
choose $g \in \EE$ such that
$g(\cl_\EE(U^i)) \subseteq W^i$ for all $i < m$.  Let $h = f - g$.
Then $h(X) \cap B(0; 1) \ne \emptyset$ (because of $x_0$),
so $\Re(h(X) \cap (-1,1)) \ne \emptyset$.
Let $z = g(p) = g(q)$.  Then $h(p) = 2-z$ and $h(q) = -2 -z$, 
so $\Re(h(X))$ meets either  $[2, \infty)$ or $(-\infty, -2]$.
But also, for each $x$, there is an $i$ such that $x \in U^i$;
so that $f(x), g(x) \in W^i$ and $h(x) \in W^i - W^i$; hence,
$\pm b \notin \Re(h(X))$.  Thus, $\Re(h(X))$ is disconnected,
so $\AA$ has a non-trivial idempotent by Lemma \ref{lemma-disc}.
\end{proof}

\begin{proofof}{Lemma \ref{lemma-R-scattered-lots}}
Immediate by Lemmas \ref{lemma-nice-sh} and
\ref{lemma-sep-LOTS-nice}.
\end{proofof}

\section{Supports of Measures}
\label{sec-supt}
We use the notion of ``support'' in two different ways to 
reduce the problem of the CSWP for a ``big'' space $X$ to
that of a ``small'' subspace.
First, we can apply it to measures which annihilate a subspace
of $C(X)$:

\begin{definition}
If $\mu$ is a Borel measure on $X$ and $\EE \subseteq C(X)$,
then $\mu \perp\EE$ means that 
$\int f \, d \mu = 0$ for all $f \in \EE$.
\end{definition}

The proof in \cite{KU} of Lemma \ref{lemma-supt-cswp}
starts with $\AA \leac C(X)$ and $\AA \ne C(X)$,
and uses a measure  $\mu\perp\AA$ to conclude that
$\supt(\mu)$ fails to have the CSWP.

Our second application uses measures associated with 
elements of the \textit{maximal ideal space}, $\MM(\AA)$.
Elements of $\MM(\AA)$
may be viewed either as maximal ideals of $\AA$, or as homomorphisms
from $\AA$ to $\CCC$.
See \cite{GA,HOR,RUD3,RUD4}.

Recall (see \cite{GA}, Theorem 4.1, or \cite{RUD3} \S5.22)
that if $\AA \leac C(X)$ and $\varphi \in \MM(\AA)$,
then there is always a regular Borel probability
measure $\nu$ on $X$ such that 
$\varphi(f) = \int f \, d \nu$ for all $f \in \AA$.
However, $\nu$ is not uniquely determined from $\varphi$, and
we wish to make the notion of ``support'' associated directly with
$\varphi$, so we make the following definition:

\begin{definition}
If $\AA \leac C(X)$ and $\varphi \in \MM(\AA)$, then a closed
$H \subseteq X$  is a \emph{pre-support} of $\varphi$
iff $|\varphi(f)| \le \|f\|_H$ for all $f \in \AA$.
A closed $H \subseteq X$  is a \emph{support} of $\varphi$
if $H$ is a pre-support but no proper subset of $H$ is a pre-support.
\end{definition}

Applying Zorn's Lemma,

\begin{lemma}
If $H$ is a pre-support of $\varphi$, then there is a closed $K \subseteq H$
which is a support of $\varphi$.
\end{lemma}

For a pre-support $H$, $\AA \res H$ need not be closed in $C(H)$,
but using $|\varphi(f)| \le \|f\|_H$, we see that
$\varphi$ defines an element of $\MM(\overline{\AA \res H})$,
and is thus represented by a probability measure on $H$.  Thus we have

\begin{lemma}
\label{lemma-supt-meas}
If $\AA \leac C(X)$,  $\varphi \in \MM(\AA)$, and $H$
is a support of $\varphi$, then there is a regular Borel probability
measure $\nu$ with $H = \supt(\nu)$ and 
$\varphi(f) = \int f \, d \nu$ for all $f \in \AA$.
\end{lemma}

For example, if 
$X \subset \CCC$ is the closed unit disc and $\AA$ is the disc algebra,
then $\MM(\AA) = X$, so every $\varphi \in \MM(\AA)$ has a singleton
as one of its supports.  Say $\varphi$ is evaluation at $0$.
Then $\{0\}$ is a support of $\varphi$, but so is
every simple closed curve in $X$ which winds around $0$.

Next we note that $\overline{\AA \res H}$ has no
non-trivial idempotents -- in fact, no non-trivial real-valued functions
(i.e., is ``anti-symmetric'' in the terminology of \cite{HS}):

\begin{lemma}
\label{lemma-antisym}
Assume that $\AA \leac C(X)$,  $\varphi \in \MM(\AA)$, and $X$
is a support of $\varphi$.  Let $f \in \AA$ be real-valued.
Then $f$ is constant.
\end{lemma}
\begin{proof}
If not, then by re-scaling, we can assume that $f: X \to [0,1]$,
and $0,1 \in f(X)$.  Let $\nu$ be as in Lemma \ref{lemma-supt-meas}.
Then $X = \supt(\nu)$, so that $0 < \varphi(f) < 1$.  Say
$\varphi(f) = 1 - 2\varepsilon$, where $0 < \varepsilon < 1/2$.
Let $g = (1 + 2\varepsilon) f$.  Then $\varphi(g) < 1$, so
$\varphi(g^n) \to 0$ as $n \to \infty$.  But also, $g(x) > (1+ \varepsilon)$
on a set of positive measure (since $X = \supt(\nu)$),
so $\varphi(g^n) = \int g^n \, d \nu \to \infty$ as $n \to \infty$. 
\end{proof}

Next, we note that a support of $\varphi$ has no isolated points
unless it is a singleton:

\begin{lemma}
\label{lemma-perfect}
Assume that $\AA \leac C(X)$,  $\varphi \in \MM(\AA)$, and $X$
is a support of $\varphi$.  Then:
\begin{itemizz}
\item[1.] $X = \Sh(\AA)$.
\item[2.] $X$ has no isolated points unless $X$ is a singleton.
\end{itemizz}
\end{lemma}
\begin{proof}
(1) follows from the definition of ``support'' as a minimal pre-support,
since $\Sh(\AA)$ is always a pre-support.
For (2), assume that $p \in X$ is isolated.  Since $\{p\}$ is
removable (Lemma \ref{remov-pt}), $X = \Sh(\AA)$ implies
that $\AA$ must have a non-trivial idempotent, contradicting
Lemma \ref{lemma-antisym}.
\end{proof}

In view of these lemmas and our methods
in Section \ref{sec-rem} for producing non-trivial
idempotents, it is important to show that in many cases,
there is some $\varphi \in \MM(\AA)$ such that at least one
of its supports is not a singleton.
Of course, this cannot be true if $\AA = C(X)$.  
We do not know if $\AA \ne C(X)$ is sufficient for obtaining
such a $\varphi$, but 
Lemma \ref{lemma-get-phi} below is a partial result in this direction
which is strong enough for our purposes.
For this, we need the following well-known theorem of \v Silov \cite{SI}
(or, see \cite{AW, HOR}):

\begin{theorem}
\label{theorem-shilov}
Suppose that  $\AA \leac C(X)$ and
every $\varphi\in\MM(\AA)$ is a point evaluation.
Then  $\AA$ contains the characteristic
function of every clopen subset of $X$.
\end{theorem}

\begin{lemma}
\label{lemma-get-phi}
Suppose that $\AA \leac C(X)$ and
$\mu$ is a non-zero regular complex Borel measure on $X$
with $\mu\perp \AA$ and $X = \supt(\mu)$. 
Suppose that some clopen $K \subseteq X$ has the
CSWP, where $\emptyset \subsetneqq K \subsetneqq X$. 
Then some $\varphi\in\MM(\AA)$ has a support which is
not a singleton.
\end{lemma}
\begin{proof}
If some $\varphi\in\MM(\AA)$ fails to be a point evaluation,
then \textit{every} support of $\varphi$ is not a singleton.
Thus, we may assume that $\MM(\AA) = X$. 
But then $\cchi_K \in \AA$ by Theorem \ref{theorem-shilov}.
Since $K$ has the CSWP, $\AA$ must contain every continuous function
which vanishes on $X \backslash K$,
contradicting $\mu\perp \AA$ (since $|\mu|(K) \ne 0$).
\end{proof}

\begin{proofof}{Theorem \ref{thm-PR-scattered-CSWP}}
Suppose that $\AA \leac C(X)$ and $\AA \ne C(X)$.
We shall derive a contradiction.
Fix a non-zero regular Borel measure $\mu$ with $\mu \perp \AA$,
let $Y  = \supt(\mu)$, and let 
$\BB = \overline{\AA\res Y} \leac C(Y)$.  Note that
$(\mu\res Y) \perp \BB$.

Now, there is a non-empty open $U \subseteq Y$
such that $\overline{U}$ is pseudo-removable,
and hence removable, since $\overline{U}$ is the
support of a measure.  It follows that $\overline{U}$
is zero-dimensional, so there is a non-empty clopen $K \subseteq Y$
with $K \subseteq \overline{U}$.  Then $K$ has the CSWP, so we can
apply Lemma \ref{lemma-get-phi} to get a 
$\varphi\in\MM(\BB)$ with some support $F \subseteq Y$ which is
not a singleton.  Applying Lemma \ref{lemma-supt-meas}, we get $\nu$
with $F = \supt(\nu)$ and 
$\varphi(f) = \int f \, d \nu$ for all $f \in \BB$.
Let $\BB'$ be the closure of $\BB \res F$.
Applying Lemmas \ref{lemma-antisym} and \ref{lemma-perfect},
all idempotents of $\BB'$ are trivial and 
$F = \Sh(\BB')$, so we have a contradiction
by Lemma \ref{lemma-reduce-PR-scattered}.
\end{proofof}

\section{Remarks and Examples}
\label{sec-conc}
We do not know whether the CSWP for $X$ and $Y$
implies the CSWP for $X\times Y$, or even for the disjoint sum of
$X$ and $Y$.

The notion of ``nice'' is closed under disjoint sums
(by Corollary \ref{cor-sum-nice}), but not under products.
For example, let $X$ be a compact separable LOTS which 
is not second countable (e.g., the double arrow space),
and let $Y$ be scattered and uncountable.  
Then $Y$ is trivially nice and $X$ is nice by Lemma
\ref{lemma-sep-LOTS-nice}, but $X \times Y$ and $X \times X$ are
not nice by:

\begin{proposition}
\label{prop-bad-product}
If $X$ is not scattered and
$Y$ is not second countable then $X\times Y$ is not nice.
\end{proposition}
\begin{proof}
For $f \in C(X \times Y)$, define $\widehat f : Y \to C(X)$ by 
$(\widehat f(y))(x) = f(x,y)$.  Then $\widehat f$ is continuous, so
$\widehat f (Y) \subseteq C(X)$ is a compact metric space, and
hence second countable.

Now, suppose that $\EE \subseteq C(X \times Y)$ is countable.
Then $\{ \widehat f  : f \in \EE \} $  is a countable family
of maps into second countable spaces.  Since $Y$ is not second
countable, there are $b,c \in Y$ with $b \ne c$ and
$\widehat h (b) = \widehat h (c)$ for all $h \in \EE$.
Thus, if $\sim$ is $\sim_\EE$, we have $(x,b) \sim (x,c)$ for
all $x\in X$.

If $\EE$ were nice, then since each $[(x,y)]$ is closed and scattered,
while $X$ is not scattered, $X \times \{b\}$ must meet uncountably many 
equivalence classes.  But then any 
$f\in C(X\times Y)$ such that $f(x,b) \ne f(x,c)$ for all $x\in X$
would contradict clause (2) of Definition \ref{def-nice}.
\end{proof}

It is easy to see from Lemma \ref{lemma-R-scat-remov}
that any finite union of removable spaces
is removable, but we do not know about products.
We also do not know if removable is equivalent to CSWP
plus totally disconnected.

If $X$ is any compact space, there are trivially compact $H$
such that the disjoint sum $X \oplus H$ fails the CSWP;
for example, $H$ itself can fail the CSWP.  Proposition
\ref{prop-compact-copy} below shows that we can always find $H$
so that $X \oplus H$ fails the NTIP.

\begin{definition}
If $f\in C(X\times Y)$, let $f_x(y) = f(x,y)$ \lparen for $x \in X$\rparen
and $f^y(x) = f(x,y)$ \lparen for $y \in Y$\rparen.
\end{definition}

\begin{lemma}
\label{lemma-ntip-prod}
Let $X$ be an arbitrary compact space.
Assume that $Y$ is compact and that $\AA \leac C(Y)$, satisfying
\begin{itemizz}
\item[1.] $\Re(h(Y))$ is connected for all $h\in\AA$.
\item[2.] Some $\varphi\in\MM(\AA)$ is not a point evaluation.
\end{itemizz}
Then there is a $\BB \leac C(X\times Y)$ such that
\begin{itemizz}
\item[a.] $\Re(f(X \times Y))$ is connected for all $f\in\BB$.
\item[b.] For each $y\in Y$, $\{f^y : f \in \BB\} = C(X)$.
\end{itemizz}
In particular $X\times Y$ fails the NTIP for all compact $X$.
\end{lemma}
\begin{proof} Let
\[
\BB = \{f \in C(X\times Y) : \exists z \in \CCC\, \forall x\in X
[ f_x \in \AA \ \&\ \varphi(f_x) =z ] \}\ \ .
\]
Observe that $\BB$ is a closed subalgebra of $C(X \times Y)$ and
$\BB$ contains all constant functions.
To prove (b), fix $y_0 \in Y$.
Since $\varphi$ is not a point evaluation, fix $h_0 \in \AA$
with $c := \varphi(h_0) \ne d := h_0(y_0)$.
Let $h = (h_0 - c)/(d - c)$;
then $h \in \AA$ and $\varphi(h) = 0$ and $h(y_0) = 1$.
Then whenever $k \in C(X)$, the function $f(x,y) = k(x)h(y)$
is in $\BB$, and $f^{y_0} = k$.

To prove that $\BB$ separates points (so that $\BB \leac C(X\times Y))$,
fix $(x_1,y_1) \ne (x_2,y_2)$.
If $y_1 \ne y_2$, fix $h \in \AA$ such that $h(y_1) \ne h(y_2)$,
and define $f(x,y) = h(y)$; then $f\in \BB$ and separates 
$(x_1,y_1) ,  (x_2,y_2)$.  
If $y_1 = y_2$, then $\BB$ separates $(x_1,y_1) ,  (x_2,y_2)$ by (b).

To prove (a), let $f\in \BB$ be an idempotent; it is enough
to show that $f$ is trivial.  By (1), each 
$f_x$ is either identically $0$ or identically $1$.
But $\varphi(f_x)$ is independent of $x$, so
$f$ itself is either identically $0$ or identically $1$.
\end{proof}

We remark that if $Y$ is connected then (1) is trivial,
whereas if $Y$ is not connected, then (1) implies (2) 
by \v Silov's Theorem \ref{theorem-shilov}.

\begin{proposition}
\label{prop-compact-copy}
If $X$ is compact, then there is a compact $Z$ containing
a clopen copy of $X$ and a $\BB \leac C(Z)$ such that
\begin{itemizz}
\item[1.] All idempotents of $\BB$ are trivial.
\item[2.] $\BB \res X = C(X)$.
\end{itemizz}
\end{proposition}
\begin{proof}
This follows from Lemma \ref{lemma-ntip-prod},
if we choose $Y$ to have an isolated point.  For example,
$Y$ can be the Cantor set plus one point
(using \cite{HS, RUD1}), or $Y$ can
be $\TTT \cup \{0\}$, where
$\AA$ is the disc algebra.
\end{proof}

We next give an example of a
PR-scattered space which is not PR:

\begin{example}
\label{ex-PR-lots}
There is a compact $X$ with a countable set $I$ of isolated points
such that if $L = X \backslash I$, then $L$ is an infinite
compact connected LOTS which does not contain a Cantor subset.
This $X$ is the support of a measure $\mu$, and $X$
PR-scattered but not PR.
\end{example}
\begin{proof}
Once $X$ is constructed, it is the support of any $\mu$
which gives positive measure to the points in $I$.
$X$ is PR-scattered because $L$ is PR
(by Lemma \ref{lemma-PR-scattered-lots}).
$X$ is not removable because it is not totally disconnected
(see Lemma \ref{lemma-remove-basic}).
Then, $X$ is also not PR because $X = \supt(\mu)$.

Let $I$ be any countably infinite set.
First,we describe $L$:  As usual, for $y,z \subseteq I$, define
$y \subseteq^* z$ iff $y \backslash z$ is finite, and define
$y =^* z$ iff $y \Delta z$ is finite.
Let $C$ be a family of subsets of $I$ satisfying:

\vbox{    
\begin{itemizz}
\item[1.] $C$ contains $\emptyset$ and $I$.
\item[2.] No two distinct elements of $C$ are $=^*$.
\item[3.] $C$ is totally ordered by $\subseteq^*$.
\item[4.] $C$ is maximal with respect to (1)(2)(3).
\end{itemizz}
}

Observe that $C$, as ordered by $\subseteq^*$, is a dense total
order with first element $\emptyset$ and last element $I$.
Let $L$ be the Dedekind completion of $C$.  Then $L$ is a compact connected
LOTS and no element of $L$
has cofinality $\omega$ from both sides, so that $L$ does not
contain a Cantor subset.
Let $X = C \cup I$.  As a subbase for the topology
of $X$, take all sets of the form $\{i\}$ and $X \backslash \{i\}$ for
$i \in I$, together with all sets of the form
$x \cup [\emptyset, x)$ and 
$(I\backslash x) \cup (x, I]$ for $x\in C$,
where $[\emptyset, x)$ and $(x,  I]$ denote intervals in $L$.
\end{proof}

Next, we show that the notion of ``nice'' provides something
not obtainable just by considering ordered spaces.
Specifically, call the compact $H$ \textit{LPR} iff
it does not contain a Cantor subset, and
whenever $\mu$ is a regular Borel measure on $H$,
$\supt(\mu)$ is homeomorphic to compact separable LOTS.
So, LPR spaces are PR by 
Lemmas \ref{lemma-sep-LOTS-nice} and \ref{lemma-nice-sh}.
Then, call the compact $X$ 
\textit{LPR-scattered} iff $X$ is scattered for the class of LPR spaces.
Note that the only concrete examples we have given for spaces with the CSWP
are contained in Theorem \ref{thm-lots-scattered} and Example \ref{ex-PR-lots},
but these spaces are actually LPR-scattered.
Example \ref{ex-nice-not-ordered} provides a space $Z'$
which has the CSWP but is not LPR-scattered.
The CSWP will follow from $Z'$ being nice and not having a Cantor subset.
To make $Z'$ not LPR-scattered, we
make sure that every non-empty 
open subset of $Z'$ fails to be a LOTS, and that
$Z'$ is separable, and hence the support of 
a measure (e.g., a countable sum of point masses).

\begin{example}
\label{ex-nice-not-ordered}
There is a compact separable nice space $Z'$ which does not contain a Cantor set
such that 
$W$ is not an ordered space whenever $W$ is a non-empty open subset of $Z'$.
\end{example}
\begin{proof}
We describe a generalization of the double arrow space construction.
Assume that $X$ is compact and has no isolated points,
and let $D$ be any finite set of ``directions''.
As a set, $Z$ will be $X \times D$, and we display
elements of $Z$ as $\Dgen_x$ instead of $(x,\Dgen)$.
For the double arrow space construction, $X = [0,1]$ 
and $D = \{\Lft, \Rt\}$, so $Z$ contains  ``left and right
copies'', $\Lft_x, \Rt_x$ for each $x\in [0,1]$.
To define the topology, choose, for each $x \in X$ and $\Dgen \in D$,
open sets $U_x^\Dgen\subseteq X \backslash \{x\}$ so that
for each $x \in X$, the $U_x^\Dgen$ are pairwise disjoint
and $\{x\} \cup \bigcup_{\Dgen \in D}U_x^\Dgen$ is also open.
In the double arrow space,
$U_x^\Rt = (x,1]$ and $U_x^\Lft = [0,x)$.
Let $\pi: Z \to X$ be the natural projection.  Give $Z$ the
topology whose subbase is all sets of the form $\pi\iv(V)$ such 
that $V$ is open in $X$ together with all sets of the
form $\{\Dgen_x\} \cup \pi\iv(U_x^\Dgen)$.
Then $Z$ is Hausdorff, $\pi$ is continuous,
and $Z$ is compact by the Alexander Subbase Lemma.
$Z$ may have isolated points, but we can discard them, forming $Z'$;
then $Z'$ has no isolated points (since $X$ has none), and
$\pi(Z') = X$.
In the double arrow space, the isolated points are $\Lft_0$ and $\Rt_1$.

If we form a base by taking finite intersections from the subbase,
then every basic open set in $Z$ is of the form
$\pi\iv(V) \cup F$, where $V$ is open in $X$ and $F$ is finite.
Thus, if $S$ is dense in $X$, then $\pi\iv(S) \cap Z'$ is dense in $Z'$.
Thus, $Z'$ will be separable whenever $X$ is separable.
$Z$ itself need not be separable, since it might have uncountably
many isolated points.

If $X$ is compact metric then $Z$ is nice:
To see that, obtain the countable
$\EE \subset C(Z)$ of Definition \ref{def-nice} by composing $\pi: Z \to X$
with any countable subfamily of $C(X)$ which separates the points of $X$.
Then $[z] = \pi\iv (\pi(z))$, which is finite.
To verify condition (2) of Definition \ref{def-nice}, it is sufficient
to fix $f\in C(Z)$ and a real $c > 0$, and show that
$\BIG_\EE(Z,f , c)$ is a finite union of equivalence classes.
If not, then
$\pi(\BIG_\EE(Z,f , c)) = \{x \in X : \diam(f(\pi\iv\{x\})) \ge c\}$
is infinite.
Since $X$ is compact, we can choose distinct $x_n \in \pi(\BIG_\EE(Z,f , c))$ 
which converge to some point $y \in X$.
Passing to a subsequence, we may assume that there is a
$\Dgen \in D$ such that each $x_n \in U_y^\Dgen$.
But then, by continuity, $\diam(f(\pi\iv\{x_n\})) < c$ for all but
finitely many $n$, a contradiction.

Call our choice of the  $U_x^\Dgen$ \textit{antisymmetric} iff
one never has both $x \in U_y^\Dgen$ and $y \in U_x^\Dgen$
for any $\Dgen \in D$ and $x,y\in X$.
Observe that this is true for the double arrow space construction.
If the choice is antisymmetric, then $Z$ does not contain a Cantor set;
in fact, no uncountable $E \subseteq Z$ can be 
second countable in its relative topology.  To prove this, 
assume that $E$ is uncountable and second countable.
We may assume (shrinking $E$ if necessary)
that there is a fixed $\Dgen \in D$ such that
$E = \{\Dgen_x : x \in G\}$, where $G$ is an uncountable subset of $X$.
Since $E$ is second countable, it is metrizable, so let
$\rho: E \times E \to \RRR$ be 
a metric on $E$ which induces the topology on $E$.
Shrinking $E$ if necessary, we may assume that there is a fixed
$\varepsilon > 0$ such that
$B(x; \varepsilon) \subseteq (\{\Dgen_x\} \cup  \pi\iv(U_x^\Dgen)) \cap E$
for all $\Dgen_x \in E$.
But then, by antisymmetry,
$\rho(\Dgen_x,\Dgen_y) \ge \varepsilon$ for all distinct
$\Dgen_x,\Dgen_y \in E$, which is impossible.

Finally, we need to obtain an example where
every non-empty open $W \subseteq Z'$ fails to be a LOTS.
Since $Z'$ is separable (assuming $X$ is), $W$ will be separable
also, so it is sufficient to make sure that $W$ fails
to be hereditarily separable (HS) (since every separable
LOTS is HS by \cite{LB}).
To do this, we modify the well-known
proof that the square of the double arrow space fails to be HS
(although this square is also not nice by
Proposition \ref{prop-bad-product}).
Let $C$ be a Cantor set in the real line with the property that
any finite subset of $C$ is linearly independent over $\QQQ$.
Let $X = C \times C$.  Observe that whenever $(x_1,x_2) \in X$,
$a,b$ are distinct positive rationals, and $t$ is a non-zero real,
then $(x_1 \pm at, x_2 \pm bt) \notin X$.
Fix $a = 1.0$ and $b = 1.1$, say;
then, for $x = (x_1,x_2) \in X$, the lines through $x$ with
slopes $\pm 1.1$ partition $X\backslash\{x\}$ into four quadrants, 
north, east, south, and west of $x$.
Let $D  = \{\No,\Ea,\So,\We\}$.
If $x = (x_1, x_2) \in X$, define:
$$
\begin{array}{c}
U_x^\No = 
    \{(w_1,w_2) \in X : w_2 > x_2 \ \&\ a(w_2 - x_2) > b |w_1 - x_1| \} \\
U_x^\So = 
    \{(w_1,w_2) \in X : w_2 < x_2 \ \&\ a(x_2 - w_2) > b |w_1 - x_1| \} \\
U_x^\Ea = 
    \{(w_1,w_2) \in X : w_1 > x_1 \ \&\ b(w_1 - x_1) > a |w_2 - x_2| \} \\
U_x^\We = 
    \{(w_1,w_2) \in X : w_1 < x_1 \ \&\ b(x_1 - w_1) > a |w_2 - x_2| \} 
\end{array} \ \ .
$$
Every non-empty open $W \subseteq Z'$
contains a subset of the form $\pi\iv V \cap Z'$ for some open $V \subseteq X$.
Inside of this open set, the points of the form $\No_x$ for $x$
on a horizontal line (i.e., with fixed $x_2$) form an uncountable discrete
set; thus, $W$ is not HS.
\end{proof}

The following modification of the proof of the 
Stone-Weierstrass theorem by de Branges \cite{BRA} (or see \cite{AW})
shows how to obtain directly (1) and (2) of
Lemma \ref{lemma-reduce-PR-scattered}.
It also provides an alternate proof of
Lemmas \ref{lemma-ntip-cswp} and \ref{lemma-supt-cswp}:

\begin{proposition}
\label{prop-de-branges}
Suppose that $\AA \leac C(X)$ and $\AA \ne C(X)$.
Then there is a non-zero complex Borel measure $\mu$ on $X$ with
$\mu \perp \AA$, such that if $F = supt(\mu)$ and $\BB = \overline{\AA\res F}$,
then all real-valued functions in $\BB$ are constant.
\end{proposition}
\begin{proof}
We identify $(C(X))^*$ with the space of measures on $X$;
note that $\|\mu\| = |\mu|(X)$.
Let $K = \{\sigma \in (C(X))^* : \sigma \perp \AA \ \&\  \|\sigma\| \le 1\}$.
Then $K$ is convex and $K$ is compact in the $\mathrm{weak}^*$ topology.
Applying the Krein-Milman Theorem, let 
$\mu\in K$ be any non-zero extreme point of $K$.
\end{proof}

We do not see how to achieve all of (1)(2)(3) of
Lemma \ref{lemma-reduce-PR-scattered} directly,
avoiding the argument of Section \ref{sec-supt}.
Note that the proof of
Proposition \ref{prop-de-branges} might result in $F = X \ne \Sh(\AA)$.
For example, let $X = \TTT \cup \{0\}$ and let
$\AA$ be the disc algebra restricted to $X$;
then $\Sh(\AA) = \TTT$. 
Let $\mu = {1 \over 2} \lambda - {1 \over 2} \delta_0$,
where $\delta_0$ is the unit point mass at $0$ and 
$\lambda$ is the usual (Haar) measure on $\TTT$.
Then $\supt(\mu) = X$ and $\mu$ is an extreme point of $K$.

Proposition \ref{prop-de-branges} suggests that we define
a notion of ``weakly removable'' (WR) by strengthening
``all idempotents of $\AA$ are trivial'' to 
``all real-valued functions in $\AA$ are constant''
in Definition \ref{def-remove}.
Then WR-scattered implies WR (as in Lemma \ref{lemma-R-scat-remov}),
and WR spaces have the CSWP, but we do not know if this increased
generality is useful, as we do not have any techniques for
producing a non-trivial real-valued function other than
by producing an idempotent (as in Lemma \ref{lemma-nice-sh}).

\newpage  

\newpage  
\section{Appendix}
\label{sec-app}
Here we collect a few remarks on Section \ref{sec-conc} 
which don't seem worth putting
in the published part of the paper, since Section \ref{sec-conc}
is itself essentially an appendix.

\medskip

Note that in Proposition \ref{prop-compact-copy},
one cannot also have
$\BB \rest (Z \backslash X) = C(Z \backslash X)$, by the following:

\begin{proposition}
If $X$ is compact, $\AA \leac C(X)$, and $X = H \cup K$,
where $H,K$ are both closed in $X$ and:
\begin{itemizz}
\item[1.] $\AA \res H = C(H)$.
\item[2.] $\AA \res K$ is dense in $C(K)$.
\end{itemizz}
Then $\AA = C(X)$.
\end{proposition}
\begin{proof}
If $g \in C(K)$ and $g \res (H \cap K) \equiv 0$, define
$e(g)\in C(X)$ so that $e(g) \res K = g $ and 
$e(g) \res H \equiv 0$.
Let
$$
\II = \{g \in C(K) : g \res (H \cap K) \equiv 0 
\ \& \ e(g) \in \AA\} \ \ .
$$
$\II$ is clearly a closed subalgebra of $C(K)$.  It is also an ideal.
To verify this, fix $g \in \II$.  
It is sufficient to show that $\{k \in C(K) : k \cdot g \in \II\}$
is dense in $C(K)$, so, by (2), it is sufficient to show that
$(f \res K) \cdot g \in \II$ whenever $f \in \AA$.
But then $e( (f \res K) \cdot g) = f \cdot e(g) \in \AA$.

Since $\II$ is a closed ideal in $C(K)$, there is a closed 
$L \subseteq K$ such that $\II = \{f \in C(K) : f \res L \equiv 0\}$.
Then $L \supseteq H \cap K$.

If $L = H \cap K$, then $\AA = C(X)$.  To verify this,
fix $f\in C(X)$.  By (1), fix $f' \in \AA$ such that 
$f' \res H = f\res H$.  Let $g = (f - f') \res K$.  Then 
$g \res (H \cap K) \equiv 0$, so $g \in \II$, so
$f - f' = e(g) \in \AA$.  Thus $f \in \AA$.

Now, suppose $L \supsetneqq H \cap K$, and fix $p \in L \setminus (H \cap K)$.
For any $f_1, f_2 \in \AA$,
if $f_1 \res H = f_2 \res H$, then $(f_1 - f_2) \res K \in \II$, so
$f_1(p) = f_2(p)$.
One can then define $\varphi : C(H) \to \CCC$ so that
$\varphi(h) = f(p)$ for some (any) $f\in \AA$ such that $f \res H = h$.
Then $\varphi \in \MM(C(H))$, so fix $q \in H$ such that
$\varphi(h) = f(q)$ for all $h \in C(H)$.
But then $p \ne q$ and  $f(p) = f(q)$ for all $f \in \AA$, a contradiction.
\end{proof}

In this proposition, one cannot weaken (1) to just
``$\AA \res H$ is dense in $C(H)$''.  For example, let $X$ be
the Cantor set and let $\AA \leac C(X)$ be the algebra defined by
Rudin \cite{RUD1}.  Then by 
Hoffman and  Singer \cite{HS}\S4, this $\AA$ is 
\textit{pervasive}; that is $\AA \res H$ is dense in $C(H)$
for \textit{every} proper closed $H \subset X$.

\medskip

We remark that for the particular $X$
obtained in Example \ref{ex-PR-lots},
one can verify the CSWP without using
the notion of PR;  $L$ has the CSWP by Theorem \ref{thm-sep-lots};
then, since $L$ is the perfect kernel
of $X$, $X$ also has the CSWP by \cite{KU}, Cor.~3.7.
However, we can modify $X$ by replacing each point of $I$
by a copy of the double arrow space; that is, we build
$Y = (X \times A) / \approx $, where $A$ is the double arrow
space and $\approx$ is the equivalence relation on $X\times A$
which identifies $\{x\} \times A$ to a point whenever $x \in L$.
Then $Y$ is also the support of a measure (since $A$ is)
and $Y$ is PR-scattered but not PR, and $Y$ is perfect.

\medskip

To verify that the $\mu$ obtained in the proof of \ref{prop-de-branges}
really satisfies the proposition,
let $\nu$ be the restriction of $\mu$ to $F = \supt(\mu)$.
Note that $\nu$ is an extreme point of
$N := \{\sigma \in (C(F))^* : \sigma \perp \BB \ \&\  \|\sigma\| \le 1\}$.
Also, $\|\nu\| = 1$ (since $\nu$ is an extreme point).

Whenever $g\in C(F)$, let $g \nu$ denote the measure defined by
$\int h \, d(g \nu) = \int h g \, d\nu$.
Note that $\|g\nu\| = \int |g| \, d|\nu|$.
Now, assume that $\BB$ contains a non-constant real-valued function $g$.
By rescaling, we may assume $g: F \to [0,1]$.  Let
$r = \|g\nu\| = \int g \, d|\nu|$; then $0 < r < 1$ 
(since $F = \supt (\nu)$), and if $s = 1-r$ and $k(x) = 1-g(x)$,
then $s = \|k \nu\|$.
But then $\nu = r\cdot (g/r )\nu + s \cdot (k/s) \nu$,
and $(g/r )\nu$ and $(k/s) \nu$ are in $N$, which is impossible,
since $\nu$ is an extreme point.

\medskip

To verify that the $\mu$ described in the paragraph following
Proposition \ref{prop-de-branges} is really an extreme point,
suppose that $\mu = (\mu^0 + \mu^1)/2$, where $\mu^0,\mu^1 \in K$.
Let $\mu^j = {1 \over 2} \nu^j - {1 \over 2} c^j \delta_0$,
where $\nu^j$ is a measure on $\TTT$ ($j = 0,1$).
Then  $\lambda = (\nu^0 + \nu^1)/2$ and $1 = (c^0 + c^1)/2$.
Note that $\|\mu^j\| = 1$ and 
$\|\mu^j\| = {1 \over 2} \| \nu^j\|  + {1 \over 2} |c^j|$.
But also, $\mu^j \perp \AA$ implies that $\int 1 \, d \mu^j = 0$,
so that $\nu^j(\TTT) = c^j$, so that $\|\nu^j\| \ge |c^j|$,
and hence $1 = \|\mu^j\|  \ge |c^j|$.  
This plus $1 = (c^0 + c^1)/2$ yields $c^0 = c^1 = 1$.
Now we have $\|\nu^j\| = 1$ and $\nu^j(\TTT) =1$, so that
$\nu^j$ is a positive measure.  But also
$f(0) = \int_\TTT f d \nu$ for all $f \in \AA$.
Since $\nu^j$ is positive, we have $\nu^j = \lambda$
(see \cite{RUD3} \S5.24), so $\mu^0 = \mu^1 = \mu$.

\end{document}